\documentclass[12pt]{article}

\usepackage{pifont}
\usepackage{geometry}
\usepackage{fleqn}
\usepackage{hyperref}
\usepackage{graphicx}
\usepackage{amsmath}
\usepackage{amsfonts}
\usepackage{amssymb}
\usepackage[T1]{fontenc}

\usepackage{latexsym}
\usepackage{times}
\usepackage{setspace}
\usepackage[bf,footnotesize]{caption}
\usepackage[displaymath]{lineno}
\usepackage{color}
\usepackage{eurosym}

\usepackage{yhmath}

\usepackage{epstopdf} 
\usepackage{epsfig}

\usepackage[T1]{fontenc}
\usepackage{authblk}

\newenvironment{proof}[1][Proof]{\noindent\textbf{#1.} }{\ \rule{0.5em}{0.5em}}

\title{Remarks on the paper "Non-existence of Shilnikov chaos in continuous-time systems"}

\author[*]{ L. Gardini\thanks{laura.gardini@uniurb.it}}
\author[**]{ D. Radi\thanks{davide.radi@unibg.it}}
\affil[*]{\emph{\small Department of Economics, Society, Politics, University of Urbino Carlo Bo}}
\affil[**]{\emph{\small Department of Management, Economics and Quantitative Methods, University of Bergamo}}

\date{}

\begin{document}
\maketitle
\begin{abstract}
The present note refers to a result proposed in \cite{ElhadjSprott2012}, and shows that the Theorem therein is not correct. We explain that a proof of that
Theorem cannot be given, as the statement is not correct, and we underline a
mistake occurring in their proof.

Since this note is supplementary to \cite{ElhadjSprott2012}, the reader should consult this paper for further explanations of the matter and the symbols used.
\end{abstract}

{\bf Keywords}: homoclinic orbit, heteroclinic orbit, chaos.

{\bf Chinese Library Classification}

{\bf 2010 Mathematics Subject Classification}

\bigskip

\section{Main section\protect\bigskip }

In this note, we show that the main result proposed in \cite{ElhadjSprott2012}, i.e. a sufficient condition for excluding the presence of
homoclinic and heteroclinic orbits, is not correct. Hence this cannot lead
to their conjecture of a fourth kind of chaos in 3D polynomial ODE systems
characterized by the non-existence of homoclinic and heteroclinic orbits.
Moreover, we remark that the conjecture can not be correct, as explained below
in detail.

\begin{figure}[h]
\begin{center}
\includegraphics[scale=1]{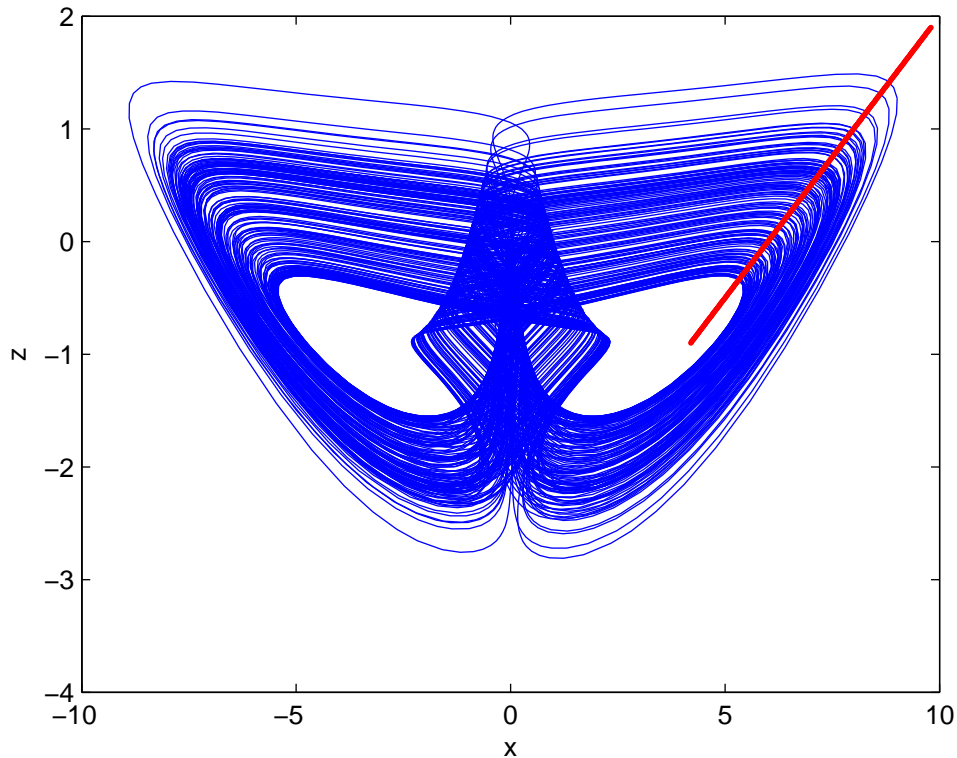}
\end{center}
\caption{Projection onto $xz$-plane of attractor obtained from system (3) in \cite{ElhadjSprott2012} for $a=40$, $b=32$ and $c=10$. The red line indicates a suitable plane for the return map.}\label{Fig1}
\end{figure}

The main result of \cite{ElhadjSprott2012} is stated in their Theorem 1, which leads to exclude the
existence of a bounded trajectory for $t<t_{0}$ in any dynamical system characterized by a vector field with at least one lower bounded component, and from their proof it
follows that this occurs independently of the existence of one or more
equilibria in the system. We notice that the proof they give also implies the \textit{non
existence of any closed orbit, i.e. limit cycle}. Moreover, \cite{ElhadjSprott2012} gives an example, in eq. (3), satisfying the assumption of the Theorem 1 and showing a chaotic
attractor illustrated in Fig. 1. From this finding, they conjecture the existence of a new type of chaos. However, given a system with a chaotic attractor (and let us assume that
this is the case shown in their example, indeed form Fig. \ref{Fig1} of this note it can be assumed that the first return map on a suitable two-dimensional surface as qualitatively shown by the red line, leads to a two-dimensional map in chaotic regime) then we have, by any definition of
chaotic system (see, e.g., \cite{ShilnikovTuraevChua12001}, \cite{ShilnikovTuraevChua22001}, \cite{GuckenheimerHolmes1983}, \cite{Kuznetsov1998} and \cite{Devaney89}), \emph{the existence of infinitely many
unstable limit cycles which densely cover the observed chaotic set.
Moreover, infinitely many homoclinic and heteroclinic orbits exist
connecting these unstable limit cycles, which also are dense in the chaotic
attractor}.
It follows that it is not possible to identify a chaotic system
characterized by the non-existence of homoclinic and heteroclinic orbits.

The paper \cite{ElhadjSprott2012} suggests that the presence of chaos is subordinated to the
existence of either homoclinic orbits \emph{of equilibria} or heteroclinic orbits
\emph{of equilibria}. Which is not correct. Indeed, the target of identifying a
chaotic system characterized by the non-existence of homoclinic and
heteroclinic orbits \textit{to one or more equilibria} may be correct.
However, in our opinion this is not interesting, as it is well known that
chaotic attractors may exist also when neither homoclinic nor heteroclinic
orbits \emph{to equilibria} are present, as it is the case in the classical Lorenz system, see, e.g. \cite{Sparrow1982}. We should emphasize that the
existence of a homoclinic orbit of an equilibrium (or heteroclinic connections
between two equilibria), under other suitable assumptions, is relevant from
a theoretical point of view, as it allows to rigorously prove the existence
of chaos and also its persistence under perturbations when the homoclinic
(heteroclinic) orbit no longer exists.

\bigskip In the following, we underline the presence of an inaccuracy in
the proof of Theorem 1 in \cite{ElhadjSprott2012}.

\begin{proof}[Remarks on the proof of Theorem 1 in \cite{ElhadjSprott2012}]
Regarding the proof of Theorem 1 in  \cite{ElhadjSprott2012}, an incorrect
part is related to the unboundedness of the trajectories. Consider a vector field 
$f=\left( f_{1},f_{2},...,f_{n}\right) ^{T}:\mathbb{R}^{n}\rightarrow 
\mathbb{R}^{n}$ belonging to class $C^{r}\left( r\geq 1\right) $, $x=\left(
x_{1},x_{2},...,x_{n}\right) ^{T}$ the state variable of the system, and $%
t\in \mathbb{R}$ the time. Assuming the existence of an $\alpha <0$, such
that for at least one $j\in \left\{ 1,2,...,n\right\} $ we have that $%
f_{j}\left( x\right) \geq \alpha $, $\forall x\in \mathbb{R}^{n}$, then (as
pointed out by the authors in eq. (2) of their paper) we have to consider the inequality%
\begin{equation}
x_{j}\left( t\right) \geq \alpha \left( t-t_{0}\right) +x_{j}\left(
t_{0}\right)   \label{cond}
\end{equation}%
\textbf{for }$\mathbf{t\geq t}_{0},$ which implies that, given a homoclinic
or heteroclinic orbit $\left( \gamma _{1}\left( t\right) ,\gamma _{2}\left(
t\right) ,...,\gamma _{n}\left( t\right) \right) $, we have $%
\lim_{t\rightarrow +\infty }\gamma _{j}\left( t\right) \geq -\infty +\gamma
_{j}\left( t_{0}\right) ,$ which is compatible with the existence of the
orbits itself. Then the authors state that the same result in \eqref{cond} holds also 
\textbf{for }$\mathbf{t<t}_{0},$ getting divergence for any orbit for $%
t\rightarrow -\infty $. However this is not correct, as for $t<t_{0}$ we
have\footnote{%
Note that, given 
\begin{equation*}
\overset{\cdot }{x(t)}_{j}=f_{j}\left( x(t)\right),
\end{equation*}
 by a simple integration from $t_{0}$ to $t$, we have
\begin{equation*}
x_{j}\left( t\right) = x_{j}\left( t_{0}\right) + \int_{t_{0}}^{t}f_{j}\left(x(s)\right)ds
\end{equation*}
Assuming $f_{j}\left( x\right)\geq\alpha$ $\forall x \in \mathbb{R}$, by basic properties of definite integrals we have that for $t>t_{0}$
\begin{equation*}
\int_{t_{0}}^{t}f_{j}\left(x(s)\right)ds \geq \alpha\left(t-t_{0}\right)
\end{equation*}
from which inequality \eqref{cond} follows, and for $t<t_{0}$
\begin{equation*}
\int_{t_{0}}^{t}f_{j}\left(x(s)\right)ds = -\int_{t}^{t_{0}}f_{j}\left(x(s)\right)ds \leq -\alpha\left(t_{0}-t\right)=\alpha\left(t-t_{0}\right)
\end{equation*}
from which we obtain inequality \eqref{RivCond}.
}%
\begin{equation}
x_{j}\left( t\right) \leq \alpha \left( t-t_{0}\right) +x_{j}\left(
t_{0}\right)   \label{RivCond}
\end{equation}%
in place of \eqref{cond}. Inequality \eqref{RivCond} implies that $%
\lim_{t\rightarrow -\infty }\gamma _{j}\left( t\right) \leq +\infty +\gamma
_{j}\left( t_{0}\right)$, which does not exclude the existence of homoclinic (or heteroclinic) orbits. 
\end{proof}

\bigskip

\section{Conclusions}

In this note we provide some arguments showing that the statement of Theorem
1 in \cite{ElhadjSprott2012} is not correct. Moreover, we show the presence of a
mistake related to the backward integration in the proof of the same Theorem.

\section*{Acknowledgments}
The authors are grateful to G.~I.~Bischi and L.~Cerboni~Baiardi  for helpful comments on this note. The usual caveats apply.

\bigskip


\end{document}